\documentclass[11pt]{amsart}
\usepackage{amsmath,amsbsy,amsfonts,epic,eepic}
\usepackage[graph,matrix]{xy}
\numberwithin{equation}{section}

\newtheorem{theorem}{Theorem}[section]  
\newtheorem{definition}[theorem]{Definition}  
\newtheorem{example}[theorem]{Example}  

\newtheorem{lemma}[theorem]{Lemma}  
\newtheorem{proposition}[theorem]{Proposition}  
\newtheorem{corollary}[theorem]{Corollary}

\def\cal{\mathcal}  

\newcommand{\F}{{\cal F}}

\newcommand{\R}{{\cal R}}

\newcommand{\Ss}{{\cal S}}

\newcommand{\flag}{\mbox{Fl}}

\newcommand{\rk}{\mbox{\rm rk}}
\newcommand{\Hom}{\mbox{\rm Hom}}
\def\kset#1{{[#1]}}
\def\nset{\kset{n}}
\def\set#1{\left\{#1\right\}}               
\def\abs#1{\left|#1\right|}                 
\def\sign{{\rm sgn}}
\def\flagify{{\F}}                          
\def\closure{{cl}}

\def\ini{{\rm In}}

\title{Annihilators of ideals of exterior algebras}
\author{Graham Denham}
\email{denham@noether.uoregon.edu}
\address{University of Oregon\\
Eugene, OR 97403}
\thanks{First author partially supported by a NSERC Postdoctoral Fellowship} 
\author{Sergey Yuzvinsky}
\email{yuz@math.uoregon.edu}
\address{University of Oregon\\
Eugene, OR 97403}
\date{\today}

\begin{document}
\begin{abstract}
The Orlik-Solomon algebra $A$ of a matroid is isomorphic to the quotient of
an exterior algebra $E$ by a defining ideal $I$.  
We find an explicit presentation of the annihilator ideal of $I$ or,
equivalently, the $E$-module dual to $A$. As an application of that
we provide 
a necessary, combinatorial condition for the algebra $A$ 
to be quadratic.  We show that this is stronger than matroid being line-closed
 thereby resolving (negatively) a conjecture by Falk. We also show that our
condition is not sufficient for the quadraticity.
\end{abstract}
\maketitle
\section{Introduction}
\bigskip
This paper concerns ideals of exterior algebras, more precisely the ideals
related to the Orlik-Solomon (OS) algebra of a matroid. 
OS algebras appeared first from theorems of Brieskorn and Orlik-Solomon as
the cohomology algebras of complements of arrangements of hyperplanes
in a complex linear space. The Orlik-Solomon theorem showed in particular that
the OS algebra of an arrangement is defined by the underlying matroid and the
definition is valid for an arbitrary 
simple matroid 
(not necessarily representable over $\mathbb C$). Such a matroid $M$ on the set
$[n]=\{1,2,\ldots n\}$ defines the exterior algebra $E$ with $n$ generators
and its graded ideal $I$ (the OS ideal of $M$). Then the OS algebra of $M$
is the graded algebra $A=E/I$.

The OS algebras have been extensively studied during the last 20 years
from different points of view and for different applications (see books
\cite{OT,OT1} and surveys \cite{Fa1,Yu}). We recall here two open problems.
The older one is to find a condition on $M$ equivalent for $A$ being
quadratic. This is important because the quadraticity of $A$ is first step to
$A$ being Koszul. The latter property is equivalent (for matroid representable
over $\mathbb C$) to
the arrangement complement being rational $K[\pi,1]$ whence it relates to other
topological properties (see \cite{PY,Fa}).
Several years ago, Falk proved that
the well-known property of $M$ being line-closed is necessary for $A$ being
quadratic and conjectured that it is also sufficient.

The other problem is a very recent one. In \cite{EPY}, the problem of
resolving $A$ as a graded $E$-module is considered.
While very little is known about free resolutions of $A$, 
the main theorem of \cite{EPY} proves that its minimal
injective resolution is linear and computes the betti numbers of $A$.
However, it leaves open the apparently hard
problem of describing this resolution explicitly.

In the present paper, we make progress on both problems by studying the
annihilator $I^0$ of $I$. 
Since $I^0$ is isomorphic to the dual module to $A$,
constructing an injective resolution of $A$ is equivalent to constructing
a free resolution of $I^0$. We explicitly find the two initial terms of this
resolution. First, we exhibit generators of $I^0$ (in fact a Gr\"obner basis).
Although we do this directly one could also use the basis of the flag space
from \cite{SV}. Then we find the generating relations among these generators.
The description of these relations is messier; although these relations
can be obtained by a deformation of the simple relations for the initial ideal
of $I^0$, this deformation needs to be chosen carefully.

To study the quadraticity of $A$ one notices that it is equivalent to the
equality 
$I^0=J^0$, where $J$ is the ideal of $E$ generated by the degree 2 component of
$I$. This equality allows us to define a property of $M$ which we call {\em
3-independent} that implies line-closed and is still necessary for the
quadraticity. This property relates to the question of 
into how many parts $[n]$ can
be broken so that each circuit has at least two points common with some of the
parts. We give an example (Example~\ref{cross}) of a line-closed matroid that
is not 3-independent, disproving Falk's conjecture. We find that 
while both ideals $I$ and $J$ are generated by pure (decomposable)
elements, their annihilators have different properties in general. The ideal
$I^0$ is also generated by pure elements (in particular the generators
mentioned above are pure); $J^0$ is not necessarily so (see Example~%
\ref{notpure}). When this happens our criterion may not work. In particular
the matroid of Example~\ref{notpure} is 3-independent but $A$ is not quadratic.

The setup of the paper is as follows. In Section 2, we study pure ideals (i.e.,
ideals generated by pure elements), in particular ideals generated by
the degree $p$ component of the OS ideal of $M$. In Section 3, we 
prove that $I^0$ is pure by exhibiting
its pure generators. We also prove that these generators form a Gr\"obner basis
and exhibit a basis of every homogeneous component of this ideal.
 In Section 4, we define $p$-independent matroids and
establish relations among the criteria of being 3-independent, line-closed, 
and quadratic. In Section
5, we give a basis of the relation space among generators of $I^0$.

\bigskip
\section{Pure ideals}
\bigskip

Let $R$ be an arbitrary commutative ring with 1,
$V$ a free module over $R$ of finite rank, and
$E=\Lambda V$ the exterior algebra of $V$. 
We view $E$ as a graded $R$-algebra with the
standard grading $E_p=\Lambda^p V$. In particular $E_1=V$.

Recall that an element $a\in E$ is called {\em pure} (decomposable) 
if it is the product of elements of degree 1. Unless $a\in V$,
the linear factors of $a$ are not 
uniquely defined. We call an ideal $I$ of $E$ {\em pure}
if it is generated by a set of pure elements.

The following classes of pure ideals constitute the main characters 
of the paper.
Let $M$ be a simple (no loops or nontrivial parallel classes) matroid 
on the set 
$[n]=\{1,2,\ldots,n\}$. Let $V$ be the $R$-module with a basis
$\{e_1,e_2,\ldots,e_n\}$. Then there is the natural homogeneous basis of $E$
consisted of all the monomials $e_S=e_{i_1}e_{i_2}\cdots e_{i_p}$ one for each 
subset $S=\{i_1,\ldots,i_p\}\subset[n]$ (in sections 1-4, we always assume that
$i_1<i_2<\cdots<i_p$). This basis allows one to define
the structure of a differential graded algebra on $E$. The differential
$\partial:E\to E$ of degree -1 is defined via $\partial e_i=1$ for every
$i=1,2,\ldots,n$ and the Leibniz rule 

$$\partial(ab)=(\partial a)b+(-1)^{\deg a}\partial b$$
for every homogeneous $a$ and every $b\in E$.
In this notation, the Orlik-Solomon (OS) ideal of $M$ is the ideal $I(M)$ of
$E$ generated by $\partial e_S$ for every {\em dependent} set $S$ of
$M$. Clearly $I(M)$ is homogeneous, $I(M)=\oplus_{p=2}^n I_p(M)$ where
$I_p(M)=I(M)\cap E_p$. Denote by $J(p,M)$ the ideal of $E$ generated by $I_r(M)$ 
with $r\leq
p$. Notice that $J(p,M)\subset J(p+1,M)$,
$J(p,M)=0$ for $p<2$, and $J({\ell},M)=I(M)$ where $\ell$
 is the rank of $M$.
A straightforward check 
shows that for every $S=\{i_1,\ldots,i_p\}\subset[n]$ we have
\begin{equation}
\label{partialeqn}
\partial 
e_S=(e_{i_2}-e_{i_1})(e_{i_3}-e_{i_1})\cdots(e_{i_{p}}-e_{i_1}).
\end{equation}
Thus all the ideals $J(p,M)$ are pure ideals.

\medskip
To each ideal $I$ of $E$ we relate its annihilator $I^0$ that is the ideal
$I^0=\{a\in E|ab=0,\ {\rm for\ every\ }b\in I\}$. Clearly $I^0$ is homogeneous if
$I$ is. There is a natural isomorphism of graded $E$-modules
$\Hom_E(E/I,E)\cong I^0$ assigning to $\phi:E/I\to E$ the element $\phi(1)$. 
Also fixing a basis in $V=E_1$ identifies $E_n$ with $R_n$ that is the trivial
$E$-module whose grading is
concentrated in degree $n$. Then the product in $E$ defines for every
$p=1,2,\ldots,n$ the nondegenerate $R$-bilinear pairing $E_p\times E_{n-p}\to R_n$. 
These
pairings define the isomorphism of $E$-modules $(E/I)^*=\Hom_R(E/I,R_n)\cong I^0$
that preserves the degree.
In particular we have the canonical isomorphism of free $R$-modules
 $I^0_p\cong(E/I)_{n-p}$. This implies also that
$(I^0)^0=I$.

Let us fix a simple matroid $M$ on $[n]$ of an arbitrary rank 
$\ell$ and put $I=I(M)$.
Extending Definition 2.64 from \cite{OT} we say
that a partition $\pi$ of $[n]$ is {\em $p$-independent} for some $p$, $3\leq
p\leq\ell+1$,
if any choice of one element from each of
 not more than $p$ arbitrary elements of $\pi$ 
forms an independent set of
$M$. It is clear that $(\ell+1)$-independent partition can have at most $\ell$
elements. Following \cite{OT} we call 
$(\ell+1)$-independent partition {\em independent}.

For any partition $\pi$ of $[n]$,
call its parts $A_1, A_2,\ldots,A_k$ (in an arbitrary order). 
 Let $\sigma\in\Sigma_n$ be the
permutation (the shuffle of $A_1$ through $A_k$) 
that puts the elements of $A_i$ before those of $A_{i+1}$
for each $i\leq k$ preserving the order inside each $A_i$ induced from $[n]$. 
Then using
this order on $A_i$ put
\begin{equation}\label{eq:defofz}
z(\pi)={\rm sign}
(\sigma)\partial(e_{A_1})\partial(e_{A_2})\cdots\partial(e_{A_k}).
\end{equation}
Clearly $z(\pi)\in E_{n-|\pi|}$.

\begin{lemma}
\label{partition}
If a partition $\pi$ of $[n]$ is $p$-independent then 
$$z(\pi)\in J(p-1,M)^0_{n-|\pi|}.$$
\end{lemma}
\proof
Let $S\subset[n]$ be a dependent set of $M$ and $|S|\leq p$. 
It suffices to prove that $z=z(\pi)$
annihilates $\partial e_S$. Since $\pi$ is $p$-independent there exists at least one
$A\in\pi$ such that $A\cap S$ contains at least two elements, say it contains 
$r$ and $s$. Then both $z$ and $\partial e_S$ are divisible by $e_r-e_s$ whence
$z\partial e_S=0$. This completes the proof.                \qed

\medskip
For $p<\ell+1$ the elements $z(\pi)$ for $p$-independent partitions do not
necessarily generate $J(p-1,M)^0$.
Moreover this ideal may not have any pure elements in its minimal degree
component, in particular it may not be a pure ideal. 
One of the celebrated Hilbert--Cohn-Vossen matroids, namely $(9_3)_2$
(see \cite{HCV}),
provides such an example when $\ell=p=3$: see Example \ref{notpure}.
The results below show that the OS ideals of matroids are special in this
respect.

\medskip
\section{OS-ideals} 
In this section, we focus our attention on the OS ideal $I=J(\ell,M)=I(M)$.
Our goal is to exhibit a Gr\"obner basis of this ideal and a generating set
consisting of certain $z(\pi)$.

It turns out that it suffices to consider only independent partitions of the
following special type.
 Let $F=(X_0=\emptyset\subset
X_1\subset X_2\subset\cdots
\subset X_{\ell}=[n])$ be a maximal flag of flats of $M$ (or a maximal chain
in the lattice $L$ of all the flats of $M$). Put
$S_i=X_i\setminus X_{i-1}$ ($i=1,2,\ldots,\ell$) and 
denote by $\pi(F)$ the ordered partition $\{S_1,\dots,S_{\ell}\}$.

\begin{lemma}
\label{independent}
For every maximal flag $F$ the partition $\pi(F)$ is independent.
\end{lemma}
\proof
In the above notation for $F$ let $i_j\in S_j$ for every $j$, $1\leq j\leq\ell$.
It suffices to prove that the set $T=\{i_1,\ldots,i_{\ell}\}$ is independent. 
If $T$ is dependent and $T_p=\{i_1,\ldots,i_p\}$ is its inclusion minimal
dependent subset
then $T_{p}\subset X_{p-1}$,
 which contradicts the
choice of $i_{p}$. The contradiction completes the proof.              \qed

\medskip
To simplify the notation put $z(F)=z(\pi(F))$ for every maximal flag $F$.
The previous two lemmas imply that $z(F)\in I^0$. Now we reduce the class of the
partitions further. Recall that a {\em broken circuit} of $M$ is a circuit
(i.e. an inclusion minimal dependent set) with its smallest element deleted.
Then a set $T
\subset [n]$ is called {\bf nbc} if there is no broken circuits lying in it.
We denote the collection of all {\bf nbc}-sets of cardinality $p$ by ${\bf
nbc}_p$.
Clearly an {\bf nbc}-set is independent.
If $T=\{i_1,\ldots,i_p\}$ (recall that $i_1<i_2<\cdots<i_p$) then it is {\bf nbc}
if and only if $i_r=\min cl(\{i_r,\ldots,i_p\})$ for each $r$, $1\leq r\leq p$.
Here by $cl(S)$ we mean the closure of $S$, i.e., the minimal flat of $M$ 
containing $S$.

It is well known that the set of monomials $e_T$ where $T$ is running through 
${\bf nbc}_p$ projects to a basis of $(E/I)_p$ under the
natural projection $E\to E/I$.

Now fix an ordered {\bf nbc}-set $T=(i_1,i_2,\cdots,i_{\ell})$
of (maximal) length $\ell$ and define the
maximal flag $\flagify(T)=
(X_0=\emptyset\subset X_1\subset\cdots\subset X_{\ell}=[n])$
via $X_p=cl(\{i_{\ell-p+1},\ldots,i_{\ell}\})$. To simplify the notation
put $z(T)=z(\flagify(T))$.

\begin{lemma}
\label{flags}
Consider the lexicographic ordering on the subsets of a cardinality $p$ of 
$[n]$ and the respective ordering on the monomials in $E_p$. Then for each 
{\bf nbc}-set
of length $\ell$ the largest ({\em leading}) monomial
of $z(T)$ is $e_{\overline T}$ where $\overline{T}=[n]\setminus T$.
\end{lemma}
\proof
Put $F=\flagify(T)$ and denote the respective partition of $[n]$ by $\pi$. 
For each
$A\in\pi$ denote by $\nu (A)$ the minimal element of $A$. Then, by definition
$(\ref{eq:defofz})$ of $z(T)$, the leading monomial of $z(T)$ is the
product of the leading monomials of $\partial(e_A)$
for each $A\in\pi$. This equals $e_{A\setminus\{\nu(A)\}}$.
 Since $T$ is {\bf nbc} we have $T=\{\nu(A)|A\in\pi\}$ 
which completes the proof.                    \qed

\medskip
Put $Z=\{z(T)|T\in{\bf nbc}_{\ell}$. 
The above lemmas imply the main result of this section
\begin{theorem}
\label{groebner}
The set $Z$ is a Gr\"obner basis of $I^0$ (with respect to the degree-lexicographic
monomial ordering).
\end{theorem}
\proof
For every set $A\subset E$, let $\ini(A)$ denote the set of all leading (i.e.
maximal in the degree-lexicographic ordering) monomials
of elements of $A$. To prove the statement we need to show that $\ini(I^0)=
\left<\ini(Z)\right>$ where the right hand side is the set of all 
monomials divisible by some
monomials from $\ini(Z)$. 
Since $Z\subset I^0$ it suffices to prove the inclusion
\begin{equation}
\label{include}
\ini(I^0)\subset \left<\ini(Z)\right>,
\end{equation}
the opposite inclusion being obvious. It is clear also that 
\begin{equation}
\label{include1}
\ini(I^0)\subset \overline{\ini(I)}
\end{equation}
where $\overline{\ini(I)}$ consists of all monomials annihilating all elements of
$\ini(I)$.
Recall that $\ini(I)$ 
consists of all the monomials divisible by some
monomials $e_U$ corresponding to broken circuits $U$.
 Thus $\overline{\ini(I)}$ consists of
monomials $e_T$ corresponding to the sets $T$
transversal to all the
broken circuits, i.e., 
intersecting each one of them nontrivially. We will call these sets
{\bf tbc}. It is easy to see that the collection of {\bf tbc}-sets
consists of the complements to {\bf nbc}-sets. In particular the inclusion
minimal  {\bf tbc}-sets have $n-\ell$ elements. Thus
\begin{equation}
\label{tbc}
\overline{\ini(I)}=\left<\ini(Z)\right>.
\end{equation}

Now \eqref{include1} and \eqref{tbc} imply \eqref{include},
which concludes the proof.                \qed

\medskip
The following corollary is a routine application of Gr\"obner basis theory.
\begin{corollary}
\label{OS-annihilator}
The set $Z$ is a generating set of $I^0$.
In particular $I^0$ is a pure ideal.
\end{corollary}

\medskip
Corollary \ref{OS-annihilator} implies a result for the OS ideal itself that 
may
have an independent use. For each {\bf nbc}-set $T$ of length $\ell$, denote by
$I_T$ the linear ideal (i.e., generated in degree 1) generated by the factors of
$z(T)$ (for instance, by $(e_i-e_{\nu(A)})$ for all $A\in \pi(T)$ 
and $i\in A\setminus\{\nu(A)\}$.
It is a particular case of a well-known property of exterior algebras that 
$(I_T)^0=Ez(T)$, where the right hand side
is the principal ideal generated by $z(T)$.

\begin{corollary}
\label{OS-ideal}
The OS ideal $I$ is the intersection of linear ideals. More precisely
$$I=\bigcap_{T\in{\bf nbc}_{\ell}} I_T.$$
\end{corollary}
\proof
Corollary \ref{OS-annihilator} gives $I^0=\sum_TEz(T)$. Thus we have
$$I=(I^0)^0=(\sum_TEz(T))^0=\bigcap_T(Ez(T))^0=\bigcap_TI_T.$$
\qed

\medskip
Theorem \ref{groebner} allows us to exhibit a basis for each $(I^0)_p$
(as a free $R$-module). For that we need to define elements $z(T)\in I^0$ for
$T\in{\bf nbc}_p$, $1\leq p\leq\ell$. The
flag $\flagify(T)=(X_0\subset X_1\subset\cdots\subset X_p)$ is defined 
exactly as
for $p=\ell$ but it is no longer maximal. More precisely, its maximal element
$X_p$ is not the maximal element of $L$. The sets $S_i=X_i\setminus X_{i-1}$
($i=1,\ldots,p$) form a partition $\pi=\pi(\flagify(T))$ of $X_p$ that is still
independent. Extend the definition of $z$ from \eqref{eq:defofz} by setting
\begin{equation}
\label{z_p}
z(T)=z(\pi)={\rm sign}(\sigma)\partial(e_{S_1})\partial(e_{S_2})\cdots
\partial(e_{S_p})e_{[n]\setminus X_p}
\end{equation}
where $\sigma$ is the shuffle of $X_p$ that puts $S_1$ through $S_p$ 
in the order of their indices and $[n]\setminus X_p$ after $S_p$.
Notice that the definition coincides with \eqref{eq:defofz} if $p=\ell$.

The properties of $z(T)$ for $|T|=\ell$ can be generalized to an arbitrary
$p$.
\begin{lemma}
\label{flags'}
For each {\bf nbc}-set $T$, the leading monomial of $z(T)$ is 
$e_{\overline{T}}$.
\end{lemma}
The proof is similar to proof of Lemma \ref{flags}.

\medskip
\begin{lemma}
\label{independent'}
For each $T\in{\bf nbc}_p$ we have
$$z(T)\in (I^0)_{n-p}.$$
\end{lemma}
\proof
Extend $\flagify(T)$ to a maximal flag $\tilde F=(X_0\subset\cdots X_p\subset
X_{p+1}\subset\cdots\subset X_{\ell})$.  By Lemmas \ref{independent}
and \ref{partition}, $z(\tilde F)\in I^0$. A straightforward check shows that
$z(T)=\pm z(\tilde F)e_{\nu(X_{p+1}\setminus X_p)}\cdots 
e_{\nu(X_{\ell}\setminus X_{\ell-1})}$. Thus $z(T)\in
I^0$ which completes the proof. \qed

\medskip
Put $Z_p=\{z(T)|T\in{\bf nbc}_p\}$.
\begin{proposition}
\label{basis}
For every $p$, the set $Z_p$ form a basis of $(I^0)_{n-p}$.
\end{proposition}
\proof
Since the $R$-module $(I^0)_{n-p}$ is isomorphic (noncanonically)
to $(E/I)_p$, we see that
the cardinality of $Z_p$ is correct. Thus it suffices to prove only that $Z_p$
generates $(I^0)_{n-p}$.

Let $e\in (I^0)_{n-p}\setminus\{0\}$. By Theorem \ref{groebner},
$\ini(e)=\pm e_Se_{\overline T}$ for $T\in{\bf nbc}_{\ell}$.
Since the set $U=[n]\setminus (S\cup\overline{T})$ 
is again {\bf nbc}, we can consider $z(U)$. Since $\ini(z(U))=e_U$, we have
$e+cz_U\in(I^0)_{n-p}$ for some $c\in R$ and $\ini(e+cz_U)<\ini(e)$.
Now the proof can be completed by induction on $\ini(e)$. \qed

\medskip
\noindent {\bf Remarks}
\begin{enumerate}
\item
In fact, the proof of Proposition \ref{basis} contains the proof that
$Z_p$ is a Gr\"obner basis of the ideal $\bigcup_{k\geq n-p}(I^0)_k$.
\item 
The signs of $z(T)$ in the formula \eqref{z_p} is chosen so that
$$e_Tz(T)=1$$ for every {\bf nbc}-set $T$.
\item
It is not hard to prove that the correspondence $F\mapsto z(F)$ induces a
map
$\flag_p\to (I^0)_{n-p}$ where $\flag$ is the graded flag module introduced by
Schechtman and Varchenko in \cite{SV} (cf. Lemma \ref{flagreln} below). 
Then Proposition~\ref{basis} can be deduced
from the results of \cite{SV}.
\end{enumerate}

\bigskip
\section{Quadratic OS-algebras}
\bigskip
As in the previous section we fix a simple
matroid $M$ on the set \hfill\break $[n]=\{1,2,\ldots,n\}$.
Its Orlik-Solomon algebra (OS algebra) is $A=A(M)=E/I(M)$. This algebra 
receives
the grading from $E$. In this section we study relations of annihilators of ideals
of $E$ with conditions for the quadraticity of $A$. The latter property means that
$I(M)$ is generated
in degree 2 (equivalently $J(2,M)=I(M)$). 

A condition sufficient for the quadraticity of $A$ was obtained by Falk in
\cite{Fa}. We recall it here. A subset $S\subset[n]$ is called {\em line-%
closed} (lcl) if it is closed with respect to 3-circuits (i.e., dependent 
sets of 3
elements). This means that if for a 3-circuit $C$ we have $|S\cap C|\geq 2$ 
then 
$C\subset S$. Clearly every flat of $M$ is lcl but the converse is false in
general. If every lcl set is a flat then $M$ itself
is called {\em line-closed} (lcl). Notice that if the rank of $M$ is 3 then 
$M$ is
lcl if and only if the line-closure of any subset $S\subset[n]$ of rank 3 is 
the whole $[n]$.

\begin{theorem}[Falk~\cite{Fa}]
\label{falk}
If $A(M)$ is quadratic then $M$ is lcl.
\end{theorem}
An alternative proof of this theorem is given below.

\medskip
Falk has conjectured that 
 the converse of Theorem \ref{falk} is also true 
and the question 
was open for a
while. The ideal annihilators allow us to solve this question negatively.

\begin{definition}
\label{indep}
Let $3\leq p\leq\ell$. The matroid $M$ is {\em $p$-independent} if any
$p$-independent partition of $[n]$ is independent.
\end{definition}

\medskip
If $\rk M=p=3$ then the following rephrasing of the above definition is
convenient. $M$ is 3-independent if  and only if any graph on the set of points
of $M$, having for each line of $M$ at least one edge lying on this line, has at
most 3 components.
The next theorem shows that being 3-independent is also necessary condition for the
quadraticity.

\begin{theorem}
\label{condition}
If $A(M)$ is quadratic then $M$ is 3-independent.
\end{theorem}
\proof
The condition of the theorem means that $J(2,M)=I(M)$ whence $(J(2,M))^0=I^0$.
Let $\pi$ be a 3-independent partition of $[n]$. By Lemma \ref{partition}
$z(\pi)\in (J(2,M))^0$ whence $z(\pi)\in I^0$. 

Assume that $\pi$ is not independent, i.e., there exists a dependent
subset $S\subset[n]$ such that $|S\cap P|\leq 1$ for every $P\in\pi$. To find a
contradiction it suffices to show that $z(\pi)\partial e_S\not=0$. Indeed 
we have equality (up to a sign)
$$z(\pi)\partial e_S=\pm z(\bar\pi)$$
where $\bar\pi$ is the partition obtained from $\pi$ by gluing together all the
elements $P\in\pi$ such that $S\cap P\not=\emptyset$. Since for 
every partition $\rho$
we have $z(\rho)\not=0$ we obtain the contradiction implying the result.
\qed

\medskip
As the next result we show that Falk's condition is implied by ours.
\begin{theorem}
\label{stronger}
If $M$ is 3-independent than it is line-closed.
\end{theorem}
\proof
Suppose $M$ is not lcl. To prove the theorem it suffices to exhibit a 3-independent
partition having more than $\ell$ elements (whence not independent). 

Let $S$ be a lcl subset of $[n]$ that is not a flat of $M$. Let $T$ be the
closure of $S$. Denote by $r$ the common rank of $S$ and $T$ and include $T$ into a
maximal flag $F$ of flats of $M$ between $T$ and $[n]$. Fixing notation, 
$F=(X_r=T\subset X_{r+1}\subset\cdots\subset X_{\ell}=[n])$. 
Since rank of $S$ is $r$ we can find a sequence $(S_0=\emptyset\subset
S_1\subset\cdots\subset S_r=S)$ of $S$ such that rank of $S_i=i$. 
Now we can define the partition $\pi$ that consists of the following
$\ell+1$ elements: 
$S_p\setminus S_{p-1}$ for $p=1,2,\ldots,r$, $S'=T\setminus S$, and $X_i\setminus
X_{i-1}$ for $i=r+1,\ldots,\ell$. 

We claim that $\pi$ is 3-independent. 
Let $X$ consist of 3 elements, one from each of 3 distinct elements of $\pi$. 
If at least one of these elements does not belong to $T$
then the independence of $X$ is obvious (cf. Lemma
\ref{partition}). The same is true if  
at least two of the elements belong to $S_{r-1}$. Thus the only interesting case is
when $X=\{a,b,c\}$ where $a\in S_{r-1},b\in Y_r=S\setminus S_{r-1},$ and $c\in
S'=T\setminus S$. In this case, we have $|X\cap S|=2$ and if $X$ is linearly 
dependent this contradicts the condition that $S$ is line-closed. Thus $\pi$ is
3-independent that concludes the proof.
           \qed

\medskip
Notice that Theorems \ref{condition} and \ref{stronger} immediately imply 
Theorem~\ref{falk}. 

The following example shows that the converse to Theorem \ref{stronger} is false,
i.e., the independence condition is strictly stronger than the line-closure 
one.

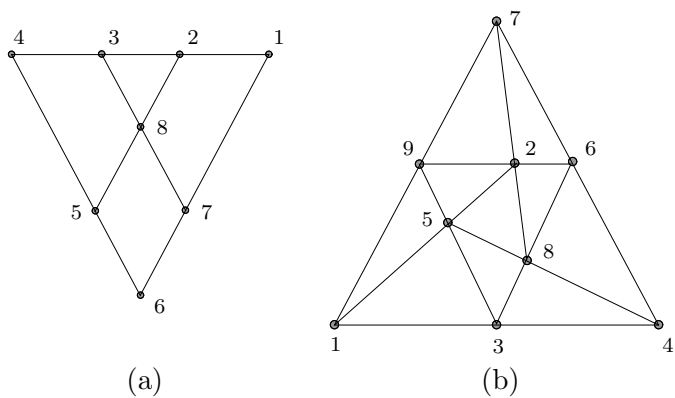
\begin{figure}
$$
\begin{array}{cc}
\raise15pt\hbox{\setlength{\unitlength}{0.00083333in}
\begingroup\makeatletter\ifx\SetFigFont\undefined%
\gdef\SetFigFont#1#2#3#4#5{%
  \reset@font\fontsize{#1}{#2pt}%
  \fontfamily{#3}\fontseries{#4}\fontshape{#5}%
  \selectfont}%
\fi\endgroup%
{\renewcommand{\dashlinestretch}{30}
\begin{picture}(1723,1767)(0,-10)
\put(41,1677){\makebox(0,0)[lb]{\smash{{{\SetFigFont{8}{9.6}{\rmdefault}{\mddefault}{\updefault}4}}}}}
\put(1124,1677){\makebox(0,0)[lb]{\smash{{{\SetFigFont{8}{9.6}{\rmdefault}{\mddefault}{\updefault}2}}}}}
\put(631,1677){\makebox(0,0)[lb]{\smash{{{\SetFigFont{8}{9.6}{\rmdefault}{\mddefault}{\updefault}3}}}}}
\put(1662,1677){\makebox(0,0)[lb]{\smash{{{\SetFigFont{8}{9.6}{\rmdefault}{\mddefault}{\updefault}1}}}}}
\put(397,588){\makebox(0,0)[lb]{\smash{{{\SetFigFont{8}{9.6}{\rmdefault}{\mddefault}{\updefault}5}}}}}
\put(1210,588){\makebox(0,0)[lb]{\smash{{{\SetFigFont{8}{9.6}{\rmdefault}{\mddefault}{\updefault}7}}}}}
\texture{55888888 88555555 5522a222 a2555555 55888888 88555555 552a2a2a 2a555555 
	55888888 88555555 55a222a2 22555555 55888888 88555555 552a2a2a 2a555555 
	55888888 88555555 5522a222 a2555555 55888888 88555555 552a2a2a 2a555555 
	55888888 88555555 55a222a2 22555555 55888888 88555555 552a2a2a 2a555555 }
\put(832,112){\shade\ellipse{42}{42}}
\put(832,112){\ellipse{42}{42}}
\put(29,1616){\shade\ellipse{42}{42}}
\put(29,1616){\ellipse{42}{42}}
\put(1633,1616){\shade\ellipse{42}{42}}
\put(1633,1616){\ellipse{42}{42}}
\put(591,1618){\shade\ellipse{42}{42}}
\put(591,1618){\ellipse{42}{42}}
\put(1077,1617){\shade\ellipse{42}{42}}
\put(1077,1617){\ellipse{42}{42}}
\put(550,639){\shade\ellipse{42}{42}}
\put(550,639){\ellipse{42}{42}}
\put(1112,642){\shade\ellipse{42}{42}}
\put(1112,642){\ellipse{42}{42}}
\put(833,1163){\shade\ellipse{42}{42}}
\put(833,1163){\ellipse{42}{42}}
\path(29,1616)(832,112)(1633,1616)(29,1616)
\path(1079,1614)(548,636)
\path(593,1614)(1115,642)
\put(918,0){\makebox(0,0)[lb]{\smash{{{\SetFigFont{8}{9.6}{\rmdefault}{\mddefault}{\updefault}6}}}}}
\put(933,1121){\makebox(0,0)[lb]{\smash{{{\SetFigFont{8}{9.6}{\rmdefault}{\mddefault}{\updefault}8}}}}}
\end{picture}
}} & \setlength{\unitlength}{0.00083333in}
\begingroup\makeatletter\ifx\SetFigFont\undefined%
\gdef\SetFigFont#1#2#3#4#5{%
  \reset@font\fontsize{#1}{#2pt}%
  \fontfamily{#3}\fontseries{#4}\fontshape{#5}%
  \selectfont}%
\fi\endgroup%
{\renewcommand{\dashlinestretch}{30}
\begin{picture}(2141,2147)(0,-10)
\texture{55888888 88555555 5522a222 a2555555 55888888 88555555 552a2a2a 2a555555 
	55888888 88555555 55a222a2 22555555 55888888 88555555 552a2a2a 2a555555 
	55888888 88555555 5522a222 a2555555 55888888 88555555 552a2a2a 2a555555 
	55888888 88555555 55a222a2 22555555 55888888 88555555 552a2a2a 2a555555 }
\put(1046,2070){\shade\ellipse{52}{52}}
\put(1046,2070){\ellipse{52}{52}}
\put(34,174){\shade\ellipse{52}{52}}
\put(34,174){\ellipse{52}{52}}
\put(2056,174){\shade\ellipse{52}{52}}
\put(2056,174){\ellipse{52}{52}}
\put(564,1179){\shade\ellipse{52}{52}}
\put(564,1179){\ellipse{52}{52}}
\put(1158,1185){\shade\ellipse{52}{52}}
\put(1158,1185){\ellipse{52}{52}}
\put(742,812){\shade\ellipse{52}{52}}
\put(742,812){\ellipse{52}{52}}
\put(1519,1193){\shade\ellipse{52}{52}}
\put(1519,1193){\ellipse{52}{52}}
\put(1046,174){\shade\ellipse{52}{52}}
\put(1046,174){\ellipse{52}{52}}
\put(1235,576){\shade\ellipse{52}{52}}
\put(1235,576){\ellipse{52}{52}}
\path(34,174)(1046,2070)(2056,174)(34,174)
\path(576,1179)(1511,1179)(1046,182)(564,1179)
\path(34,174)(1172,1185)
\path(1046,2070)(1234,578)
\path(2056,174)(753,807)
\put(1223,1231){\makebox(0,0)[lb]{\smash{{{\SetFigFont{8}{9.6}{\rmdefault}{\mddefault}{\updefault}2}}}}}
\put(1597,1231){\makebox(0,0)[lb]{\smash{{{\SetFigFont{8}{9.6}{\rmdefault}{\mddefault}{\updefault}6}}}}}
\put(579,792){\makebox(0,0)[lb]{\smash{{{\SetFigFont{8}{9.6}{\rmdefault}{\mddefault}{\updefault}5}}}}}
\put(2079,4){\makebox(0,0)[lb]{\smash{{{\SetFigFont{8}{9.6}{\rmdefault}{\mddefault}{\updefault}4}}}}}
\put(1022,0){\makebox(0,0)[lb]{\smash{{{\SetFigFont{8}{9.6}{\rmdefault}{\mddefault}{\updefault}3}}}}}
\put(8,8){\makebox(0,0)[lb]{\smash{{{\SetFigFont{8}{9.6}{\rmdefault}{\mddefault}{\updefault}1}}}}}
\put(1126,2047){\makebox(0,0)[lb]{\smash{{{\SetFigFont{8}{9.6}{\rmdefault}{\mddefault}{\updefault}7}}}}}
\put(462,1235){\makebox(0,0)[lb]{\smash{{{\SetFigFont{8}{9.6}{\rmdefault}{\mddefault}{\updefault}9}}}}}
\put(1335,591){\makebox(0,0)[lb]{\smash{{{\SetFigFont{8}{9.6}{\rmdefault}{\mddefault}{\updefault}8}}}}}
\end{picture}
} \\
\hbox{(a)} & \hbox{(b)}\\
\end{array}
$$
\caption{(a) $M$ of Example~\ref{cross}; (b) the $(9_3)_2$ matroid;}
\label{932}
\end{figure}

\begin{example}
\label{cross}
Let $n=8$ and $M$ the matroid of rank 3 shown in Figure~\ref{932}(a).
Its underlying set is $[8]=\{1,\ldots,8\}$, and it has the following
maximal dependent sets of rank 2: 
$$\{\{1,2,3,4\}, \{1,6,7\}, \{2,5,8\}, \{3,7,8\}, \{4,5,6\}\}.$$
By directly checking each independent set of rank $3$,
it is easy to see that the line-closure of each of them is $[8]$; that is, $M$
is line-closed. On the other hand, consider the partition 
$$\pi=\{\{1,3,7\}, \{2,4,5\}, \{6\}, \{8\}\}.$$
 Since none of the dependent 3-sets  intersect
non-trivially with more than two elements of $\pi$, it is 3-independent. 
However, since $\pi$
has $4$ elements, it is not independent. Hence $M$ is not 3-independent.
Notice that the graph consisting of $2$ triangles with the 
vertices $1,3,7$ and 
$2,4,5$ plus two isolated points $6$ and $8$ has an edge on each line of the
matroid.
\end{example}
The next example shows in turn that the converse to Theorem \ref{condition} is
also false: if $M$ is $3$-independent, $A(M)$ need not to be quadratic.

\begin{example}\label{notpure}
Consider the rank $3$ matroid $M$ with nine points shown in 
Figure~\ref{932}(b). 
This matroid first appeared in the book \cite{HCV} by Hilbert and Cohn-Vossen,
where it is called the $(9_3)_2$ configuration.
Since $M$ is fixed we put $J=J(2,M)$ for the length of this example.

That $A(M)$ is not quadratic can be seen by the following simple computation
(that applies in fact to the other two $9_3$ configurations from \cite{HCV}).
The ideal $J$ is generated by 9 elements from $E_2$. Since each of them is
annihilated by a 2-dimensional subspace of $E_1$ we have
\begin{equation}
\label{LHS}
\dim J_3\leq7\times 9=63.
\end{equation}
On the other hand, the Hilbert series of $A(M)$ should be divisible by $1+t$ whence
a simple computation gives
$H(A(M),t)=1+9t+27t^2+19t^3$. This implies
\begin{equation}
\label{RHS}
\dim I_3={9\choose 3}-19=65.
\end{equation}
Thus $J\not=I$.

Proving that $M$ is 3-independent requires more work. Let us say for
convenience that a subset $T$ of the vertices of $M$ {\em represent}
a 3-circuit $S$
if $|T\cap S|\geq 2$. A useful observation about $M$ is that it has an
automorphism
of order 9 that can be written as a permutation $\tau=[261594837]$ in 
two-line notation.  In particular, it is transitive on the vertices.
Now we state the following for an 
arbitrary subset $T$ of the vertices of $M$ (it suffices to consider
cases where
$1\in T$).

\begin{enumerate}
\item[(1)] If $\abs{T}=3$ then it represents at least one circuit.
\end{enumerate}

Indeed the only vertices that do not represent a circuit with 1 are 6 and 8 and
$\{6,8\}$ represents a circuit.

\begin{enumerate}
\item[(2)]
If $\abs{T}=4$ then it represents at least 2 and at most 5 circuits. 
\end{enumerate}

There are two cases. Suppose $T$ contains a circuit $S$ (one can assume 
$S=\{1,3,4\}$ or $S=\{1,2,5\}$). Then by inspection any other vertex represents
another circuit with one of the vertices of $S$. The first statement follows.
Clearly in this case $T$ can represent at most 4 circuits.

Now suppose $T$ does not contain a circuit. Since two distinct subsets of $T$
of cardinality $3$ cannot represent the same circuit, again the first
assertion must hold.  
On the other hand, $M$ viewed as a graph does not have a complete
subgraph of with $4$ vertices, and this implies the second assertion.

\begin{enumerate}
\item[(3)]
If $\abs{T}=5$ then it represents at most $7$ circuits. If $\abs{T}=6$ 
then this number is $8$.
\end{enumerate}

This follows from (1) and (the first assertion of) (2) by passing to
the complements. 

Now suppose $\pi$ is a $3$-independent partition of the set of vertices of $M$.
This means that the elements of $\pi$ represent each of the $9$ circuits.
Considering all the possibilities for the cardinalities of elements of $\pi$
if $\abs{\pi}=4$, i.e. $\{6,1,1,1\},\ \{5,2,1,1\},\ \{4,3,1,1\},\ \{4,2,2,1\}$,
we deduce from (the second statement of) (2) and (3) that the elements of $\pi$
can represent at most $8$ elements. Thus $\abs{\pi}\leq 3$ whence $\pi$ is
independent. This shows that $M$ is $3$-independent.

In fact using Macaulay~2 \cite{GS} we can obtain a basis of $(J^0)_5$ and show
that $(J^0)_5$ does not contain any pure elements (whence $J^0$ is not a pure
ideal). To describe this basis let
$$
x=e_1e_2e_3e_6e_8-e_1e_2e_4e_6e_9+e_1e_2e_4e_8e_9-e_1
e_2e_6e_8e_9-e_2e_3e_5e_6e_9,$$
$p=(1-\tau)(1+\tau^3+\tau^6)x$, and $q=\tau p$.
Then $p$ and $q$ form a basis for $(J^0)_5$.

For any element $r\in E_5$, let $\Ss(r)$ be the support of $r$: that is,
the set of all subsets
$\set{i_1,i_2,\ldots,i_5}$ that index nonzero monomials $e_{i_1}\cdots e_{i_5}$
in $r$.  Suppose $r$ is pure, that is $r=a_1a_2\cdots a_5$ where each
$a_i\in E_1$.  Then the coefficients of the monomials that make up $r$
are the minors of a $9\times5$ matrix whose $i$th column is given by $a_i$.
In particular, then, the set $\Ss(r)$ is the set of bases of the matroid
on the set $\kset{9}$ given by the rows of this matrix.

Consider a nonzero linear combination $r=\alpha p+\beta q$.  We show that
$r$ cannot be pure as follows.  We can assume that $\alpha\neq0$, by 
replacing $p$ and $q$ with $\tau p$ and $\tau q$, respectively, otherwise.
Then put
$B=\set{1,2,6,8,9}$.
 By writing out $r$ we find 
$B_1\in\Ss(r)$.  
However, neither $B_1\cup\set{5}\setminus\set{9}$ nor
$B_1\cup\set{7}\setminus\set{9}$ are in $\Ss(r)$.  Thus $\Ss(r)$ does not
satisfy the main axiom for the set of bases of a matroid.  It follows that
$(J^0)_5$ contains no pure elements.
\end{example}

\section{A Presentation of $I^0$}

The goal of this section is to exhibit a presentation of $I^0$ as an $E$-module.
Using the generating set $Z$ we only need to describe basic relations among
elements of $Z$.

First we recall a little about the minimal resolution of the
initial ideal $J=In(I^0)$. For that it is convenient to start with the monomial
ideal $J'$ of the polynomial ring $\R=R[x_1,\ldots,x_n]$ generated by
the same monomials as $J$ (i.e., $e_{\bar T}$ where $T\in{\bf nbc}_{\ell}$)
 where $x_i$ is substituted for $e_i$ ($i=1,\ldots,n$). 
The ideal $J'$ is the Stanley-Reisner ideal
of a simplicial complex $\Delta'$ that is the canonical Alexander dual to the {\bf
nbc}-complex $\Delta$ (cf., for example, \cite{EPY} and \cite{RW}). 
By the latter we mwan the complex on $[n]$ whose simplexes are the {\bf nbc}-sets.
Keeping this in mind, an explicit realization 
of the  minimal resolution of $J'$ (as an $\R$-module) can be 
given in terms of the homology of links of $\Delta$ or the local homology of
 the lattice of the least common multiples of {\bf tbc}-sets.
These constructions can be deduced from \cite{Yu1} or 
\cite{RW}. Then a minimal resolution of $J$ as an $E$-module can be constructed
using \cite{AAH} or \cite{EPY}.

For a presentation of $J$
 we do not need any of those general constructions.
A generating set of $J$ is mentioned above.
In order to describe basic relations among the generators
denote by ${\bf nbc}'$ the subset of ${\bf nbc}_{\ell-1}$ consisting
 of sets $S$ lying
each in at least two ${\bf nbc}$-bases. Notice that ${\bf nbc}'$
can be identified with the rank 2 part of the lattice of the least common
multiples of {\bf tbc}-sets of cardinality $n-\ell$.
Then put for every  $S\in{\bf nbc}'$
$$N(S)=\{i\in[n]\setminus S|S\cup\{i\}\in{\bf nbc}_{\ell}\}$$
and denote by $i(S)$ the minimal element of $N(S)$. 
The natural basic relations among these generators can be
 divided in two kinds. Relations
 $r(S,i)$ of the first kind are 
parametrized by all the pairs $(S,i)$ with $S$ as
above and $i\in N(S)\setminus\{i(S)\}$. 
The relation $r(S,i)$ is
\begin{equation}
\label{d_1}
\epsilon([n]\setminus S,i)
e_ie_{{\bar S\setminus\{i\}}}-e_{i(S)}e_{{\bar S\setminus\{i_S\}}}=0
\end{equation}
where $\epsilon (T,i)=|\{j\in T|j<i\}|$.
The relations $s_{T,j}$ of the second kind are
parametrized by pairs $(T,j)$ where $T\in{\bf nbc}_{\ell}$,
$j\in T$, and the relation $s(T,j)$ is 
$e_je_{\bar T}=0$.

There is a couple of features 
of the relations of the first kind that we will use later. 
First, they are partitioned in parts parametrized by ${\bf nbc}'$ and the part
corresponding to $S$ contains 
  $|N(S)|-1$
relations. Second, instead of the particular relations \eqref{d_1} that
we have chosen as basic, one can take any similar relations 
using a set of pairs 
$\{\{i,j\}|i,j\in N(S)\}$ as long the graph on $N(S)$ 
whose edges are these pairs is a tree
containing all the elements of $N(S)$.

The minimal resolution of $J=In(I^0)$ gives some information about the minimal
resolution of $I^0$. As it was observed in \cite{EPY}, $J$ is a flat
degeneration 
of $I^0$. It also was proved in \cite{EPY} that the minimal resolutions of $J$
and $I^0$ are linear
whence these two resolutions have the same
dimensions of the corresponding terms. 

We already know
that generating sets of $J$ and $I^0$ are both
in one-to-one correspondences with ${\bf nbc}_{\ell}$.

The situation is harder though
with linear relations (of the first kind) among $z(T)$.
In general there are no sufficiently many relations among them 
involving only pairs (similar to \eqref{d_1}). Because of that we 
describe the basic relations in two stages.  
On the first stage,
we establish linear relations between pairs of $z(F)$ that
are indexed by flags more general than $\flagify(T)$ for $T\in{\bf nbc}_{\ell}$. 
 We are able to find enough of them so they
form a basis of the whole space of relations.
However, these $z(F)$ are not
in general linearly independent over $R$, and the second stage is occupied with
expressing them in terms of the standard $z(T)$ 
that form  an $R$-basis of $(I^0)_{\ell}$.

We call two maximal flags {\em close} if and only if they differ by at most one
flat. 

Let $S$ and $T$ be disjoint subsets of $\nset$, and set $s=\min S$, $t=\min T$.
Let $\sigma$ be the shuffle of $S$ and $T$.
Using \eqref{partialeqn} we have
$$
\partial(e_S)(e_s-e_t)\partial(e_T)=\sign(\sigma)\partial(e_{S\cup T}).
$$
It follows that if $\pi$ is a partition of $[n]$ with
$\pi=\set{A_1,A_2,\ldots,A_k}$, $s=\min A_i$, and $t=\min A_{i+1}$, then
\begin{equation}
\label{merge}
(e_s-e_t)z(\pi)=(-1)^{\sum_{j\leq i}(\abs{A_j}-1)}z(\tilde{\pi}),
\end{equation}
where
$\tilde{\pi}$ 
is the partition with $k-1$ parts, obtained from $\pi$ by
joining $A_i$ with $A_{i+1}$.  

Suppose that $F$ and $F'$ are distinct, close, maximal flags: that is,
$F=(X_0,X_1,\ldots,X_\ell)$ and $F'$ differs from $F$ by replacing $X_k$
with $X'_k$ for a single $k$, $0<k<\ell$.  Let $\tilde{\pi}$ be the partition
obtained from $\pi(F)$ or (equivalently) $\pi(F')$ by joining
the $k$th and $k+1$st parts.  We find that $z(F)$ and $z(F')$ are both
divisors of $z(\tilde{\pi})$ using \eqref{merge}.  In order to be
more specific, let $a=\min X_k\setminus X_{k-1}$ and 
$b=\min X_{k+1}\setminus X_k$.  Obtain $a'$ and $b'$ similarly using $X'_k$
instead of $X_k$.  Then
\begin{lemma}
\label{linrelations}
In the notation above, 
\begin{equation}\label{eq:linrelns}
(e_a-e_b)z(F)-(-1)^{\abs{X_k}-\abs{X'_k}}(e_{a'}-e_{b'})z(F')=0.
\end{equation}
\end{lemma}
\proof
Immediate from equation \eqref{merge}.
\qed

Now we define the class of flags and monomials we will use.
Recall that if the elements of the $\ell$-tuple
$U=(u_1,u_2,\ldots,u_\ell)$ comprise a base, 
 $\flagify(U)$
denotes the flag whose flat in rank $p$ is the closure of the last $p$ 
elements of $U$.  Previously we considered $U$ which were {\bf nbc}-bases
written in increasing order.  Call such a $U$ and its flag $\flagify(U)$
{\em standard}.  In what follows, we will make use of 
other orderings of {\bf nbc}-bases.

Any maximal flag
$F$ defines an ordered base, which we denote by $\varphi(F)$,
the $\ell$-tuple whose 
$p$th element is equal to $\min(X_p\setminus X_{p-1})$, for $1\leq p\leq\ell$.
Clearly $\varphi(F)$ is an ordered base with $F=\flagify(\varphi(F))$.
Recall from Lemma~\ref{flags} that
if $U$ is standard, then $\varphi(\flagify(U))=U$.
If $U$ is not in increasing order, however,
$\varphi(\flagify(U))$ may not equal $U$.  
Using the similar definition from \cite{BV},
call any ordered base $U$ {\em neat} if $\varphi(\flagify(U))=U$.

\begin{lemma}
\label{charofneat}
A neat ordered base is {\bf nbc} (but not necessarily increasing).
\end{lemma}
\proof
If $U$ is increasing and neat, then $U$ is an ordered {\bf nbc}-base.
If $U=(u_1,u_2,\ldots,u_\ell)$ is neat but $u_{i-1}>u_i$ for some
$i$, then we show that the ordered base
$U'$ obtained by transposing $u=u_{i-1}$ and $v=u_i$
is also neat.  From this the conclusion will follow by induction.

Let $X=\closure\set{u_{i+1},\ldots,u_\ell}$.  Since $U$ is neat,
$v=\min(X\vee\set{v}\setminus X)$ and 
$u=\min(X\vee\set{v,u}\setminus X\vee\set{v})$.
The inclusion
$$
X\vee\set{u}\setminus X\subset \left(X\vee\set{v,u}\right)
\setminus \left(X\vee\set{v}\right)
$$
implies $u=\min X\vee\set{u}\setminus X$
and
\begin{eqnarray*}
\min(X\vee\set{v,u}\setminus X\vee\set{u})&\geq&
\min(X\vee\set{v,u}\setminus X)\\
&=&v.
\end{eqnarray*}
Since $v$ is in the set on the left, we have equality,
and $U'$ is neat.
\qed

\medskip
Now fix $S\in{\bf nbc}'$.
We shall construct a tree whose vertices are labeled
by $N(S)$, and whose edges give $\abs{N(S)}-1$ linear relations 
of the type \eqref{eq:linrelns}.

Write $S=\set{s_1,s_2,\ldots,s_{\ell-1}}$ with $s_1<\cdots<s_{\ell-1}$.
Some notation is needed for ordered {\bf nbc}-bases obtained by adding
an element of $N(S)$ to $S$: for $a\in N(S)$ and $0\leq k\leq\ell$, 
put $S(a,k)=(s_1,\ldots,s_{k-1},a,s_k,\ldots,s_{\ell-1})$.   If $a>s_{k-1}$,
that is, if $a$ appears at or to the left of its place in increasing order,
call $S(a,k)$ {\em early}.

In what follows, set $X^k=\closure\set{s_k,s_{k+1},\ldots,s_{\ell-1}}$,
for $1\leq k\leq\ell-1$.  We define a graph $\Gamma=\Gamma(S)$ 
from $S$ and $N(S)$ as
follows.  The vertices of $\Gamma$ are those sets $S(a,k)$ which are 
neat and early.  The edges are the pairs $S(a,k)$ and
$S(b,k+1)$ for which $a\vee X^k=b\vee X^k$.  Note that, 
if two vertices of $\Gamma$ are joined by an edge, then the flags of those two
vertices are close.

\begin{example}
\label{crossII}
Let $M$ be the matroid of Example~\ref{cross} (Figure~\ref{932}(a)).  
Of the 14 {\bf nbc}-bases,
only $125$, $135$, $145$, and $157$ contain the set $S=\set{1,5}$, so
in this case $N(S)=\set{2,3,4,7}$.  One finds that the vertices of $\Gamma(S)$
are $S(2,1)$, $S(a,2)$ for each $a\in N(S)$, and $S(7,3)$.  $\Gamma$ is a tree,
rooted at $S(2,1)=(2,1,5)$, shown in Figure~\ref{trees}.  
(The ordered sets are written with the element of $N(S)$ emphasized.)

Similarly, taking $S=\set{1,3}$ gives $N(S)=\set{5,6,7,8}$, and the 
resulting tree is also shown in Figure~\ref{trees}.
\begin{figure}
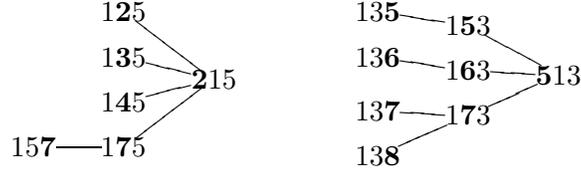

$$
\xygraph{[]~{<0.6cm,0cm>:<0cm,0cm>}
[]!{(2,-1.5)}*{{\bf 2}15}="215"-@[]!{!(-2,1.5)}*{1{\bf 2}5}="125"
[d]*{1{\bf 3}5}-@{-}"215"
"125"[dd]*{1{\bf 4}5}-@{-}"215"
"125"[]!{(-2,-3)}*{15{\bf 7}}-@{-}[]!{!(2,0)}
*{1{\bf 7}5}-@{-}"215"
}\qquad\qquad
\xygraph{[]~{<0.6cm,0cm>:<0cm,0cm>}
[]!{!(2,-0.3)}
*{1{\bf5}3}="153"-@{-}[]!{!(2,-1.1)}
*{{\bf5}13}="513"
"153"-@{-}[]!{!(-2,0.3)}*{13{\bf5}}="135"
[d]*{13{\bf6}}-@{-}[]!{!(2,-0.3)}
*{1{\bf6}3}-@{-}"513"
"135"[]!{!(0,-2.2)}*{13{\bf7}}-@{-}[]!{!(2,-0.1)}
*{1{\bf7}3}="173"-@{-}"513"
"135"[]!{!(0,-3.2)}*{13{\bf8}}-@{-}"173"
}
$$
\caption{Trees $\Gamma(S)$ for $S=\set{1,5}$ and $S=\set{1,3}$}
\label{trees}
\end{figure}
\end{example}

\begin{lemma}
\label{neatcondition}
If $a\in N(S)$, $S(a,k)$ is early, and $a=i(S)=\min N(S)\cap(X^k\vee\set{a}
\setminus X^k)$, then
$S(a,k)$ is neat.
\end{lemma}
\proof
Since $S(a,k)$ is early and $a\in N(S)$, 
$S(a,r)$ is standard for some $r\geq k$: let $F=\flagify(S(a,r))$.
The first $\ell-r$ and last $k$ flats of $S(a,k)$ 
agree with those of $F$, so we only need check the condition in ranks
$\ell-r+1$ through $\ell-k+2$.  In rank $\ell-i$, 
for $i$ with $k\leq i<r$, we have 
$s_i=\min X^i\vee\set{a}$ since $X^i\vee\set{a}$ is a flat of $F$; therefore
$s_i=\min X^i\setminus X^{i+1}$.  In the remaining rank, $\ell-k+2$, 
let $b=\min (X^k\vee\set{a})\setminus X^k$.  Then $S(b,k)$ is neat, so
by Lemma~\ref{charofneat}, $b\in N(S)$.  By hypothesis, then, $b=a$.
\qed

\medskip
As before, let $S\in{\bf nbc}'$ and
 $a=\min N(S)$.
\begin{proposition}
$\Gamma(S)$ is a tree, rooted at $S(a,1)$.
The leaves of $\Gamma(S)$ are the standard {\bf nbc}-bases indexed
by $N(S)$.
\end{proposition}
\proof
We claim that any vertex $S(b,k+1)$ for $k\geq1$ is connected to 
exactly one vertex $S(c,k)$.  Since such an $S(c,k)$ must be neat, 
the only possibility would be 
$$
c=\min(N(S)\cap(X^k\vee\set{b}\setminus X^k)).
$$
We have only to show that, for this choice of $c$, $S(c,k)$ is actually
in $\Gamma$.

$S\cup\set{b}$ is a {\bf nbc}-base; let $F=(Y_0<Y_1<\cdots<Y_\ell)$ be
its standard flag.  Since $S(b,k+1)$ is early, the last $k+1$ flats of
$F$ and $\flagify(S(b,k+1))$ agree; in particular, $X^{k-1}\vee\set{b}=
Y_{l-k+2}$.  Since $F$ is a standard flag, $s_{k-1}=\min Y_{\ell-k+2}$, so
$s_{k-1}<c$, which means $S(c,k)$ is early.

$S(c,k)$ is also neat, by Lemma~\ref{neatcondition}, 
and so (by definition) it is a vertex of $\Gamma$.  This proves the first claim.

To prove the second claim, observe that any $S(c,k)$ that is standard
is certainly in $\Gamma$.  Suppose that such a vertex $S(c,k)$
is connected to some $S(b,k+1)$ (hence is not a leaf).  
Then $c<s_k<b$.  We have
\begin{eqnarray*}
\lefteqn{
\min(X^k\vee\set{c}\setminus X^{k+1})\;=\;\min(X^k\vee\set{b}\setminus X^{k+1})
}\\
&=&\min\left(X^k\vee\set{b}\setminus X^{k+1}\vee\set{b}\right)\cup
\left(X^{k+1}\vee\set{b}\setminus X^{k+1}\right)\\
&=&\min\set{s_k,b}=s_k,
\end{eqnarray*}
since $S(b,k+1)$ is neat.  The same argument applied to $S(c,k)$ shows that the
minimum is $c$, a contradiction.

Conversely, if $S(c,k)$ is not standard, then the argument of 
Lemma~\ref{charofneat} shows that $S(c,k+1)$ is also in $\Gamma$.  
The two are
connected, so $S(c,k)$ is not a leaf.
\qed

\medskip
Given such a tree $\Gamma(S)$, 
construct a graph ${\mathfrak t}={\mathfrak t}(S)$ with vertices $N(S)$ and
an edge for a pair $(a,b)$, $a\neq b$,
if and only if $S(a,k)$ and $S(b,k+1)$ are joined by
an edge for some $k$ in $\Gamma(S)$. Then ${\mathfrak t}$ is also a tree, since its 
edges must have $a<b$.  Its construction 
can also be described algorithmically: for each non-minimal 
element $b\in N(S)$, write
the base $S\cup\set{b}$ in increasing order to give a standard ordered base
$S(b,r)$ for some $r$.  Now
move $b$ to the left to get $S(b,k)$, for $k=r-1$, $k=r-2$, $\ldots$,
as long as $S(b,k)$
remains neat.  For some greater $k$, $S(b,k)$ will not be neat.  Then
$\varphi(\flagify(S(b,k)))=S(a,k)$ for some smaller $a\in N(S)$.  Connect
$a$ to $b$, and repeat until the resulting graph is connected.

By construction, we have:
\begin{proposition}
\label{treerelns}
For each $S\in{\bf nbc}_{\ell-1}$, each edge of the graph ${\mathfrak t}(S)$
on the set $N(S)$, constructed above, gives a pair of close flags,
$F=\flagify(S(a,k))$ and $F'=\flagify(S(b,k+1))$, for some $k$.
\end{proposition}

\medskip
Then Lemma~\ref{linrelations} describes a linear relation
between $z(F)$ and $z(F')$. These are the relations of the first kind. 
It is much easier to generalize the relations of the second kind.
For each $T\in{\bf nbc}_{\ell}$, write $\pi(T)=\{A_1,\ldots,A_{\ell}\}$.
Then we have $n-\ell$ relations 
\begin{equation}
\label{trivial}
(e_j-e_{\nu(A_i)})z(T)=0
\end{equation}
where $j\in A_i\setminus\{\nu(A_i)\}$, $i=1,2,\ldots,\ell$.

What is left for us to observe is that all these relations, of the first and the
second kind together,
form an $R$-basis of the
relation space.

\begin{theorem}
\label{relbasis}
The relations \eqref{eq:linrelns}  where a pair $\{F,F'\}$
runs through all edges of all the 
graphs ${\mathfrak t}(S)$ ($S\in{\bf nbc}'$)
together with the relations \eqref{trivial} where 
$T\in{\bf nbc}_{\ell}$
 form an $R$-basis of the relation space.
\end{theorem}

\proof Observe that for each $S$ we have precisely $|N(S)|-1$
relations of the first kind. 
Since this number coincides with the number of relations of the first kind for
$J$ and the same is true for the relations of the second type, the total number
of the relations is correct. Now for each of those relations we can obtain a
relation for the generators of $J$ by taking the initial monomials of $z(F)$,
$z(F')$, $z(T)$
and the coefficients. The latter relations are $R$-linearly independent.
The nontrivial part of the reasoning for that is that the involving pairs 
for the relations of the first kind with fixed $S$ from ${\bf nbc}'$
 form a tree. Thus the 
relations for the elements from $Z$ are also independent 
which completes the proof.           \qed

\medskip

\begin{example}
Continuing the previous Example~\ref{crossII}, we obtain for $S=\set{1,5}$
and $S=\set{1,3}$, the following respective trees ${\mathfrak t}(S)$:
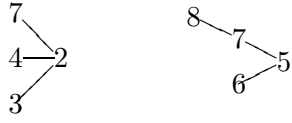
$$
\xygraph{[]~{<0cm,0.6cm>:<0cm,0cm>}
*{3}-@{-}[dr]
*{2}="2"-@{-}[u]
*{4}[r]
*{7}-@{-}"2"
}\qquad
\qquad
\raise10pt\hbox{\xygraph{[]~{<0cm,0.6cm>:<0cm,0cm>}
*{6}-@{-}[]!{!(0.5,-1)}
*{5}-@{-}[]!{!(0.5,1)}
*{7}-@{-}[]!{!(0.5,1)}
*{8}
}}
$$
For example, the edge $(5,6)$ of ${\mathfrak t}(\set{1,3})$ comes from the edge
joining $S(6,2)$ to $S(5,1)$ in 
Figure~\ref{trees}.  The corresponding flags
$F=(\emptyset<\set{3}<\set{3,6}<\kset{8})$ and $F'=(\emptyset<
\set{3}<\set{1,2,3,4}<\kset{8})$ are close.  By Lemma~\ref{linrelations},
$$(e_6-e_2)z(F)-(e_1-e_5)z(F')=0.$$
\end{example}

The remaining step of
writing a presentation of $I^0$ is to reexpress flags of elements of the
trees $\Gamma$ in terms of standard flags.  
For that it is convenient to use linear dependencies (over $R$) that
are essentially equivalent to the flag relations (2.1.1) of \cite{SV},
\begin{lemma}
\label{flagreln}
The elements $z(F)$ satisfy the equation
$$
\sum_{X_{i-1}\subset Y\subset X_{i+1}}z(X_0\subset\cdots\subset 
X_{i-1}\subset Y\subset X_{i+1}\subset\cdots\subset X_\ell)=0,
$$
for any maximal flag $(X_0\subset X_1\subset\cdots\subset X_\ell)$ 
and any $i$, $0<i<\ell$.
\end{lemma}
\proof
Let $A_1, A_2, \ldots, A_k$ denote the sets $Y\setminus X_{i-1}$, as
$Y$ ranges over all flats between $X_{i-1}$ and $X_{i+1}$.  Note that
these sets
partition $X_{i+1}\setminus X_{i-1}$, since $Y$
is the closure of $X_{i-1}\cup\set{y}$ for any $y\in Y\setminus X_{i-1}$.
Now let $\sigma$ be a permutation of $\nset$ so that the
elements of $X_{i+1}\setminus X_{i-1}$ are consecutive, and the elements
$A_j$ come before $A_{j+1}$, for all $1\leq j<k$.  

To prove the statement
it is enough to factor out those terms that appear in each summand
and show instead that
$$
\sum_{j=1}^k\sign(\sigma^j)\partial(e_{A_j})\partial(e_{A_1}e_{A_2}\cdots
\widehat{e_{A_j}}\cdots e_{A_k})=0,
$$
where $\sigma^j$ is the permutation obtained from $\sigma$ by moving the
block $A_j$ in front of $A_1$.  Put $a_j=\abs{A_j}$ for each $j$.  Since
$\partial$ satisfies the Leibniz rule, the sum becomes
\begin{eqnarray*}
\lefteqn{
\sign(\sigma)\sum_{j=1}^k(-1)^{a_j\sum_{p<j}a_p}
\partial(e_{A_j})\partial(e_{A_1}e_{A_2}\cdots
\widehat{e_{A_j}}\cdots e_{A_k})}\\
&=&\sign(\sigma)\partial\left(
\sum_{j=1}^k(-1)^{a_j\sum_{p<j}a_p}
\partial(e_{A_j})e_{A_1}e_{A_2}\cdots
\widehat{e_{A_j}}\cdots e_{A_k}\right)\\
&=&\sign(\sigma)\partial(\partial(e_{A_1}e_{A_2}\cdots e_{A_k}))\\
&=&0.
\end{eqnarray*}

\qed

\medskip

\begin{proposition}
\label{changeofbasis}
Let $S$ be an independent set with $\abs{S}=\ell-1$.
If $S(a,k)$ is in $\Gamma(S)$, then
$$
z(S(a,k))=\sum_{S(b,j)}(-1)^{j-k}z(S(b,j)),
$$
where the sum is taken over all standard $S(b,j)$ that are descendents
of $S(a,k)$ in $\Gamma(S)$.
\end{proposition}
\proof
Choose any $S(a,k)$ in $\Gamma$.  If it is not standard, then we claim that
$z(S(a,k))=-\sum_b z(S(b,k+1))$, taking 
the sum over all $S(b,k+1)$ in $\Gamma$ 
connected to $S(a,k)$.  By Lemma~\ref{flagreln}, 
$$
z(S(a,k))=-\sum_i z(F_i)
$$
where sum is taken over all flags $F_i$ that differ from $\flagify(S(a,k))$ in
rank $\ell-k$, that is, agreeing at all flats except $X^k$.  
The image of $\varphi$ is always neat, so it follows from 
Lemma~\ref{charofneat} that $\varphi(F_i)=S(b_i,k+1)$ for some $b_i\in N(S)$.
By the neatness, $b_i=\min(X^{k+1}\vee\set{b_i}\setminus X^{k+1})$ and
$a=\min(X^k\vee\set{a}\setminus X^k)$.  The first set is contained in the
second, so $b_i\geq a$.  By
assumption $S(a,k)$ is not standard: $a>s_k$.  This means
$S(b_i,k+1)$ is early, hence an element of $\Gamma$, and the claim is proven.

The proposition now follows by induction.
\qed

\medskip
By combining Proposition~\ref{treerelns} with Proposition~\ref{changeofbasis},
our presentation is described completely.

\begin{example}
Continuing the previous example, the same relation becomes 
$$
(e_6-e_2)\left(-z(136)\right)-
(e_1-e_5)\left(z(135)+z(136)+
z(137)+z(138)\right)=0
$$
in terms of the standard basis for $I^0$.
\end{example}

\end{document}